\documentclass[12pt]{article}
\usepackage{amssymb} 
\usepackage[all]{xy}

\newtheorem{thm}{Theorem}[section]
\newtheorem{cor}[thm]{Corollary}
\newtheorem{lem}[thm]{Lemma}

\newtheorem{prob}[thm]{Problem}
\newtheorem{definition}[thm]{Definition}
\newenvironment{defn}{\begin{definition}\rm}{\end{definition}}
\newtheorem{example}[thm]{Example}
\newenvironment{exmp}{\begin{example}\rm}{\end{example}}
\newtheorem{remark}[thm]{Remark}

\newcommand{\eqr}[1]{~\mbox{$(${\rm \ref{#1}}$)$}}
\newcommand{\Section}[1]{\section{#1}\setcounter{equation}{0}}

\newcommand{\cL}{{\mathcal L}}
\newcommand{\cZ}{{\mathcal Z}}
\newcommand{\so}{{\mathcal O}}
\newcommand{\tD}{\widetilde{\mathcal D}}
\newcommand{\tB}{\widetilde{B}}

\newcommand{\A}{{\mathbb A}}
\newcommand{\C}{{\mathbb C}}
\newcommand{\K}{{\mathbb K}}

\newcommand{\pn}[1]{{\mathbb P}^{#1}}
\newcommand{\mat}{M\! at_{n\times  n}}
\newcommand{\matmn}{M\! at_{m\times  n}}
\newcommand{\matmp}{M\! at_{m\times  p}}
\newcommand{\rank}{{\rm rank}\,}
\newcommand{\row}{{\rm rowsp}\,}

\newcommand{\G}{{\rm Grass}}
\newcommand{\zwei}[2]{\left[ \begin{array}{c}
                   #1 \\ #2 \end{array} \right]}
\newcommand{\vier}[4]{\left[ \begin{array}{ccc}
                   #1 &\;& #2 \\ #3 &\;& #4 \end{array} \right]}

\newcommand{\openbox}{\leavevmode
  \hbox to.77778em{%
    \hfil\vrule
  \vbox to.675em{\hrule width.6em\vfil\hrule}%
  \vrule\hfil}} \newcommand{\proofname}{Proof}
\newenvironment{proof}[1][\proofname]{\par\normalfont
  \trivlist\item[\hskip\labelsep\itshape #1:]\ignorespaces
  }{\hspace*{1cm}\hspace*{\fill}\openbox \medskip\endtrivlist}

\title{Zero Assignment, Pole Placement and \\
  Matrix Extension Problems: \\
  A Common Point of View}

\author{Meeyoung Kim\\
  Department of Mathematics\\
  Michigan State University \\
  East Lansing, MI 48824-1027 \\
  {\em e-mail:} kim@math.msu.edu %
  \and Joachim Rosenthal\thanks{Supported in part by NSF
    grant  DMS-96-10389.} \\
  Department of Mathematics,\\
  University of Notre Dame,\\ Notre
  Dame, Indiana 46556, USA.\\
  {\em e-mail:} Rosenthal.1@nd.edu%
  \and Xiaochang Alex Wang\\
  Department of Mathematics\\
  Texas Tech University\\
  Lubbock, TX 79409-1024\\
  {\em e-mail:} mdxia@ttacs.ttu.edu}

\begin{document}\maketitle

\begin{abstract}
  The paper studies a general inverse eigenvalue problem which
  contains as special cases many well studied pole placement and
  matrix extension problems.  It is shown that the studied
  problem corresponds on the geometric side to a central
  projection from some projective variety. The degree for this
  variety is computed in the critical dimension.
\end{abstract}

\Section{Introduction and motivational examples}

Let $\K$ be an arbitrary field and consider matrices of size
$E,A$ of size $n\times n$ and matrices $B,H$ of size $n\times m$.
These matrices define the discrete dynamical
\begin{equation}
  \label{dis-sys}
  Ex_{t+1}+Ax_t+Hu_{t+1}+Bu_t=0.
\end{equation}

Consider the vector space $\matmn$ consisting of all $m\times n$
matrices defined over $\K$. Let $\cL\subset \matmn$ be a linear
subspace of dimension~$d$. This paper will be devoted to the
following `constrained' pole placement question:
\begin{prob}                          \label{problem}
  Given an arbitrary monic polynomial $\varphi(s)\in\K[s]$ of
  degree $n$. Is there a feedback law of the form $u_t=Fx_t$
  having the following properties:
  \begin{enumerate}
  \item The closed loop system
   \begin{equation}
    \label{clo-sys}
    (E+HF)x_{t+1}+(A+BF)x_t=0
  \end{equation}
  has characteristic polynomial $\varphi(s)$.
\item The feedback satisfies the constraint $F\in\cL$.
  \end{enumerate}
\end{prob}

By definition the characteristic polynomial of the
system\eqr{clo-sys} is the unique monic polynomial which is a
scalar multiple of $\det \left[ s(E+HF)+(A+BF)\right]$.

Note that the set of monic polynomials of degree $n$ can be
identified with the vector space~$\K^n$.  If
Problem~\ref{problem} has a positive answer for all monic
polynomials $\varphi(s)\in\K^n$ of degree $n$, then we will say
that the system\eqr{dis-sys} is arbitrarily pole assignable in
the class of feedback compensators $\cL$.  If for a generic set
of monic polynomials $\varphi(s)\in\K^n$ of degree $n$
Problem~\ref{problem} has a positive answer, then we will say
that system\eqr{dis-sys} is generically pole assignable in the
class of feedback compensators~$\cL$.

A dimension argument immediately reveals that\eqr{dis-sys} is
generically pole assignable only if $\dim\cL\geq n$. Another
natural necessary condition for generic pole assignability is the
left primeness of the matrix pencil $[ sE\!+\! A\ \ sH\!+\!  B]$.
This last condition is satisfied for a generic set of matrices
$E,A,B,H$.

If the pencil $[ sE\!+\! A\ \ sH\!+\!  B]$ is left prime then the
transfer function $(sE+A)^{-1}(sH+B)$ defines a system of
McMillan degree $n$ which has the generic controllability indices
and every $n\times m$ transfer function of McMillan degree $n$
with the generic controllability indices is of the form
$(sE+A)^{-1}(sH+B)$.  In terms of transfer functions the problem
therefore asks: Given a $n\times m$ transfer function
$(sE+A)^{-1}(sH+B)$ of McMillan degree $n$ which has the generic
controllability indices and given a polynomial
$\varphi(s)\in\K[s]$ of degree $n$ whose roots are disjoint from
the roots of $\det(sE+A)$, is there a $m\times n$ matrix $F$ with
$F\in\cL$ and having the property that the zeroes of the transfer
function $\left[ I_n+(sE+A)^{-1}(sH+B)F\right] $ coincide with
the roots of $\varphi(s)$?

The following set of examples show that Problem~\ref{problem} is
very general indeed and it contains in special cases many well
studied pole placement and matrix extension problems.

\begin{exmp}
  Let $E=I_n$, $H=0$ and let $\cL=\matmn$. In this situation
  Problem~\ref{problem} consists of the well known state feedback
  pole placement problem. In this case Problem~\ref{problem} has
  a solution if and only if the matrices $A,B$ form a
  controllable pair. In other word the genericity condition
  $\rank\left[B, AB, \ \ldots \ A^{n-1}B\right]=n$ has to be
  satisfied. This last condition is equivalent with the left
  primeness of the matrix pencil $\left[ sI_n\!+\! A\ B\right]$.
\end{exmp}
\begin{exmp}
  Consider the static output feedback pole placement problem over
  the complex numbers $\C$:
  \begin{equation}
    \label{pole}
    \dot{x}=Ax+Bu,\ y=Cx, \ \ x\in\C^n, u\in\C^m \mbox{ and
    }y\in\C^p.
  \end{equation}
  The problem asks for a static feedback law $u=Ky$ such that the
  closed loop system
  $$
  \dot{x}=(A+BKC)x,\ y=Cx
  $$
  has some desired closed loop characteristic polynomial. One
  immediately verifies that Problem~\ref{problem} covers this
  situation if one chooses $E=I_n$, $H=0$ and $\cL:=\{ KC\mid
  K\in\matmp\}$. Over the complex numbers the main result in this
  area of research was given by Brockett and Byrnes~\cite{br81}.
  It states:
  \begin{thm}                    \label{BrBy}
    If $n\leq mp=\dim\cL$ then for a generic set of matrices
$A,B,C,$ the system\eqr{pole}
is arbitrarily pole assignable.
    Moreover if $n=mp$ then when counted with multiplicities
    there are exactly as many solutions as the degree of the
    complex Grassmannian variety $\G (m,\C^{m+p})$ once embedded via
    the Pl\"ucker embedding.
  \end{thm}
\end{exmp}

\begin{exmp}                     \label{Ex-5}
  Let $m=n$, $E=I_n$, $H=0$ and $B=I_n$. In this case
  Problem~\ref{problem} asks for conditions which guarantee that
  the characteristic map
 \begin{equation}
   \label{ext-pro}
   \chi_A:\ \cL \longrightarrow\K^n,\hspace{3mm}F \longmapsto
 \det(sI+A+F)
 \end{equation}
 is surjective or at least `generically' surjective. This general
 matrix extension problem contains itself many of the matrix
 completion problems as they were studied
 in~\cite{ba90a,fr72,fr77a,go95}.

 The main result in the situation of Example~\ref{Ex-5} has been
 derived in~\cite{he97}. It states:
  \begin{thm}
    If the base field $\K$ is algebraically closed then for a
    generic set of matrices $A\in\mat$ the characteristic
    map\eqr{ext-pro} is dominant (generically surjective) if and
    only if
    \begin{enumerate}
    \item $\dim\cL\geq n.$
    \item There must be at least one element $L\in\cL$ whose
      trace tr$(L)\neq 0$, i.e. $\cL\not\subset sl_n.$
    \end{enumerate}
  \end{thm}
\end{exmp}

The main result of this paper (Theorem~\ref{main}) will show that
over an algebraically closed field system\eqr{dis-sys} is
generically pole assignable for a generic set of matrices
$E,A,B,H$ if and only if $\dim \cL\geq n$.

The paper is structured as follows: In Section~\ref{Sec-compact}
we will introduce a natural compactification of the linear space
$\cL$ which we will denote by $\bar{\cL}$. In order to prove the
main theorem we will show that one has a characteristic map
$\chi$ defined on a Zariski open set of the variety~$\bar{\cL}$.
Geometrically $\chi$ describes a central projection from the
variety~$\bar{\cL}$ to the projective space~${\mathbb P}^n$. As a
consequence the number of solutions in the critical dimension,
i.e. in the situation where $\dim\cL=n$, is equal to
$\deg\bar{\cL}$ when counted with multiplicities and when taken
into account some possible `infinite solutions'.

The degree of the variety $\bar{\cL}$ is of crucial importance
for the understanding of the characteristic map $\chi$.  In
Section~\ref{Sec-special} we compute the degree of $\bar{\cL}$ in
many special cases. As a corollary we will rediscover several
matrix completion results as they were derived earlier
in~\cite{by93,fr72,fr77a,go95}.

Finally in Section~\ref{Sec-generic} we will compute the degree
of $\bar{\cL}$ for a generic subspace $\cL\subset\mat$.\bigskip

We would like to thank M.A. de Cataldo for useful discussions
which helped to prove the results presented in
Section~\ref{Sec-generic}.  The first author would like to thank
the Max-Planck-Institut f\"ur Mathematik in Bonn for its support
and hospitality. A large part of the work was carried out during
her stay at that institute.

\Section{Compactification of the problem}    \label{Sec-compact}

The inverse eigenvalue problem formulated in
Problem~\ref{problem} describes an intersection problem in the
linear variety $\cL$. In order to invoke results from
intersection theory~\cite{fu84} it is important to understand the
intersection at the `boundary' of $\cL$. What is needed is a good
compactification of $\cL$. It turns out that
Problem~\ref{problem} induces in a natural way a compactification
and we will explain this in the sequel.

The closed loop characteristic polynomial can be written as
\begin{eqnarray}
  \label{clo-cha}
  \det \left[ s(E+HF)+(A+BF)\right]&=&
  \det [sE\! +\! A\mid sH\! +\! B]\left[\begin{array}{c}I_n\\ F
  \end{array}\right]\nonumber\\
  &=&\det\vier{I_m}{F}{-\! sH\! -\! B}{sE\! +\! A}.
\end{eqnarray}
Following an idea introduced by Brockett and Byrnes~\cite{br81}
for the static output pole placement problem we will identify
$\row [I_m\ F]$ with an element of $\G(m,\K^{m+n})$. In this way
we have natural embeddings
$$
\cL\longrightarrow \G(m,\K^{m+n}) \longrightarrow {\mathbb
  P}\left(\wedge^{m} \K^{m+n}\right)=:\pn{N}.
$$
\begin{defn}
  Let $\bar{\cL}$ be the projective closure of $\cL$.
\end{defn}

By definition $\bar{\cL}$ is a projective variety of dimension
$\dim\bar{\cL}=\dim\cL$. The remainder of the paper will be
devoted to a large extend in the study of this variety.  In order
to have a general idea of how the projective closure of $\cL$ is
defined, we start with an illustrative example.
\begin{exmp}
  Let $\cL\in M\! at_{3\times 3}$ be defined by
\begin{equation}
\cL=\left\{ \left[\begin{array}{ccc}
      a&b&0\\
      c&a&b\\
      0&c&a
\end{array}\right]\right\}
\end{equation}
where $a,\ b,\ c\in\K$ are arbitrary elements. Then
for fixed $a,b,c$
\begin{equation}\label{grass36}
\left[\begin{array}{cccccc}
    1&0&0&a&b&0\\
    0&1&0&c&a&b\\
    0&0&1&0&c&a
\end{array}\right]
\end{equation}
is a point in ${\rm Grass} (3,\K^6)$. Let $z_{ijk}$ be the full
size minor of (\ref{grass36}) consisting of the $i$th, $j$th,
$k$th columns. Then $\{z_{ijk}\}$ are the Pl\"{u}cker coordinates
of ${\rm Grass} (3,\K^6)$ in $\pn{19}$. $\cL$ is defined by $6$
linear equations of its entries.  In terms of the Pl\"{u}cker
coordinates, they become
\begin{equation}\label{dfeqs1}
\begin{array}{lll}
z_{234}=-z_{135},&z_{234}=z_{126},&z_{235}=-z_{136},\\
z_{125}=-z_{134},&z_{124}=0,      &z_{236}=0.
\end{array}
\end{equation}
$\cL$ has $9$ minors of size $2\times 2$, but there are only $6$
monomials of degree $2$ of $a,b,c$:
$$a^2,b^2,c^2,ab,ac,bc.$$
So there are $3$ linear relations among
the $2\times 2$ minors. In terms of the Pl\"{u}cker coordinates,
they are
\begin{equation}\label{dfeqs2}
\begin{array}{lll}
z_{146}=-z_{245},&z_{345}=z_{156},&z_{346}=-z_{256}.
\end{array}
\end{equation}
The monomials $a^2,b^2,c^2,ab,ac,bc$ are not
algebraically
independent, they
satisfy the relation
\begin{equation}\label{segre}
\rank \det \left[\begin{array}{ccc}
a^2&ab &ac \\
ab &b^2&bc \\
ac &bc &c^2
\end{array}\right]\leq 1
\end{equation}
i.e. all the $2\times 2$ minors of (\ref{segre}) are zero, which
induce $6$ quadratic relations among the $2\times 2$ minors of
$\cL$:
\begin{equation}\label{dfeqs3}
\begin{array}{ll}
z_{346}^2+z_{246}z_{356}=0,&z_{146}^2+z_{246}z_{145}=0,\\
(z_{246}+z_{345})^2-z_{356}z_{145}=0,&z_{246}(z_{246}
+z_{345})-z_{346}z_{146}=0,\\
z_{346}(z_{246}+z_{345})+z_{356}z_{146}=0,
&z_{146}(z_{246}+z_{345})+z_{145}z_{346}=0.
\end{array}
\end{equation}
It is not hard to show that $\bar{\cL}$ is defined by
(\ref{dfeqs1}), (\ref{dfeqs2}), and (\ref{dfeqs3}) in ${\rm
  Grass} (3,\K^6)\subset \pn{19}$.
\end{exmp}

Note that every element in $\bar{\cL}$ can be simply represented
by a subspace of the form $\row [F_1\ F_2]$, where the $m\times
m$ matrix $F_1$ is not necessarily invertible. Row span$\, [F_1\
F_2]$ describes an element of $\cL$ if and only if $F_1$ is
invertible. Note that a characteristic equation is even defined
if $F_1$ is singular unless the polynomial in\eqr{clo-cha} is the
zero polynomial.

Let $f_i, \; i=0,\ldots,N$ be the Pl\"ucker coordinates of $\row
[F_1\ F_2]$.  In terms of the Pl\"{u}cker coordinates the
characteristic equation can then be written up to a constant
factor as:
\begin{equation}    \label{eqw2}
\det\vier{F_1}{F_2}{-\! sH\! -\! B}{sE\! +\! A}=
\sum_{i=0}^N f_ip_i(s),
\end{equation}
where the $p_i(s)$ represent up to sign and order the full size
minors of $[-\! sH\! -\! B\ \ sE\! +\! A ]$.

Let $\cZ \subset \pn{N}$ be the linear subspace defined by
\begin{equation}    \label{center}
\cZ =\{z\in {\mathbb P}^N|\sum_{i=0}^N p_i(s)z_i=0\}.
\end{equation}
Identify a closed loop characteristic polynomial $\varphi(s)$
with a point in $\pn{n}$. In analogy to the situation of the
static pole placement problem considered in~\cite{br81,wa92}
(compare also with~\cite[Section~5]{ro94}) one has a well defined
characteristic map
\begin{equation}    \label{central}
   \begin{array}{cccc}
   \chi: &\bar{\cL}-\cZ & \longrightarrow &{\mathbb P}^n \\
    & \row [F_1\ F_2] & \longmapsto &\sum_{i=0}^N f_ip_i(s).
   \end{array}
\end{equation}
It will turn out that surjectiveness of the map $\chi$ will imply
the generic pole assignability of system\eqr{dis-sys} in the
class of compensators $\cL$. The geometric properties of the map
$\chi$ are as follows:
\begin{thm}
  The map $\chi$ defines a central projection. In particular if
  $\cZ \cap \bar{\cL}=\emptyset$ and $\dim\cL=n$ then $\chi$ is
  surjective of mapping degree equal to the degree of the variety
  $\bar{\cL}$.
\end{thm}

The proof for this theorem is identical to the one given
in~\cite{wa92}. In the algebraic geometry literature (see
e.g.~\cite{ha92b,mu76}) $\chi$ is sometimes referred to as a {\em
  projection of $\cL$ from the center $\cZ$ to $\pn{n}$} and
$\cZ\cap \bar{\cL}$ is sometimes referred to as the {\em base
  locus}. Of course the interesting part of the theorem occurs
when $\cZ \cap \bar{\cL}=\emptyset$ since in this situation very
specific information on the number of solutions is provided.  If
$\cZ \cap \bar{\cL}=\emptyset$ and $\dim\cL=n$ then one says that
$\chi$ describes a {\em finite morphism} from the projective
variety $\bar{\cL}$ onto the projective space $\pn{N}$.

In analogy to the situation of the static pole placement
problem~\cite{br81,wa92} and the dynamic pole placement
problem~\cite{ro94} we introduce a definition for this important
situation:
\begin{defn}
  A particular system $E,A,B,H$ is called $\cL$-{\em
    nondegenerate} if $\cZ \cap \bar{\cL}=\emptyset$. A system
  which is not $\cL$-nondegenerate will be called
  $\cL$-degenerate.
\end{defn}
In general it will always happen that certain systems $E,A,B,H$
are $\cL$-degenerate.  The next theorem shows that if the
dimension of $\cL$ is not too large then the set of matrices
$E,A,B,H$ which are $\cL$-degenerate are contained in a proper
algebraic subset when viewed as a subset in the vector space
$\K^{2m(m+n)}$.
\begin{lem}
  Assume the base field $\K$ is algebraically closed. If $\dim\cL
  >n$ then every system $E,A,B,H$ is $\cL$-degenerate. If
  $\dim\cL \leq n$ then a generic set of systems $E,A,B,H$ is
  $\cL$-nondegenerate.
\end{lem}
\begin{proof}
  If $\dim\cL >n$ then $\cZ \cap \bar{\cL}$ is nonempty by the
  (projective) dimension theorem (see e.g.~\cite{mu76}) and the
  fact that $\dim\cZ\geq N-n-1$.

  Assume now that $\dim\cL \leq n$. Identify the set of systems
  $E,A,B,H$ with the vector space $\K^{2m(m+n)}$. In analogy to
  the proof of~\cite[Lemma 5.3]{ro94} we compute the dimension of
  the coincidence set
  \begin{equation}
    \label{coinci}
    {\mathcal S}:=\left\{ (F_1,F_2;\ E,A,B,H)\in\bar{\cL}\times
      \K^{2m(m+n)} \mid \det\vier{F_1}{F_2}{-\! sH\! -\! B}{sE\!
        +\! A}=0 \right\} .
  \end{equation}
  Using the same arguments as in~\cite{ro94} one computes
  $$
  \dim{\mathcal S}=\dim\bar{\cL}+2m(m+n)-n-1.
  $$
  Since $\bar{\cL}$ is projective the projection onto the
  second factor (namely $\K^{2m(m+n)}$) is an algebraic set by
  the main theorem of elimination theory (see e.g.~\cite{mu76}).
  This projection can result in an algebraic set of dimension at
  most $\dim{\mathcal S}$. The claim therefore follows.
\end{proof}

We are now in a position to state one of the main theorems of
this paper:
\begin{thm}      \label{main}
  Assume the base field $\K$ is algebraically closed. Let
  $\cL\subset\mat$ be a fixed subspace. Then the map $\chi$
  introduced in\eqr{central} is surjective for a generic set of
  matrices $E,A,B,H$ if and only if $\dim\cL\geq n$. If
  $\dim\cL=n$ then for a generic set of matrices $E,A,B,H$ the
  intersection $\cZ \cap \bar{\cL}=\emptyset$ and the
  characteristic map $\chi$ describes a finite morphism of
  mapping degree which is equal to the degree of the variety
  $\bar{\cL}$.
\end{thm}
\begin{proof}
  (Compare with~\cite[Theorem 2.14]{ra96}). If $\dim\cL< n$ then
  a simple dimension argument shows that $\chi$ cannot be
  surjective. We therefore will assume that $\dim\cL:=d\geq n$.

  Consider once more the coincidence set ${\mathcal S}$
  introduced in\eqr{coinci} and consider the projection onto
  $\K^{2m(m+n)}$. For a generic point inside $\K^{2m(m+n)}$ the
  fiber of the projection has dimension equal to $d-n-1$.  Let
  $E,A,B,H$ be a system whose fiber has this dimension and let
  $\cZ$ be the corresponding center as defined in\eqr{center}. By
  construction we have that
  $$
  \dim \cZ\cap \bar{\cL }=\dim \cL-n-1=d-n-1.
  $$
  In particular if $\dim\cL=n$ then $\cZ \cap
  \bar{\cL}=\emptyset$ and the characteristic map $\chi$ is a
  finite morphism.

  If $d>n$ choose a subspace $H\subset \pn{N}$ having codimension
  $d-n$ inside $\pn{N}$ and having the property that
\begin{equation}
\bar{\cL}\cap\cZ\cap H=\emptyset.
\end{equation}
Such a subspace $H$ exists by~\cite[Corollary (2.29)]{mu76}. Let
$\pi_1:\bar{\cL}\rightarrow \pn{d}$ be the central projection
with center $\cZ\cap H$ and let
$\pi_2:\pn{d}-\pi_1(\cZ)\rightarrow \pn{n}$ be the central
projection with center $\pi_1(\cZ)$.  Then $\pi_1$ is a finite
morphism which is surjective over $\K$. $\pi_2$ is a linear map,
it is surjective as well and
$$
\chi=\pi_2\circ \pi_1.
$$
It follows that $\chi$ is surjective as soon as $\dim\cL\geq
n$.
\end{proof}
Theorem~\ref{main} assumes that the field is algebraically
closed. For general fields it is often possible to deduce some
results by considering the corresponding question over the
algebraic closure. The following results is of this sort:
\begin{cor}
  If the degree of the variety $\bar{\cL}$ defined over the
  complex numbers $\C$  is odd and if $\dim\cL\geq n$ then $\chi$
  is also surjective over the real numbers ${\mathbb R}$ for a
  generic set of real matrices $E,A,B,H$.
\end{cor}
\begin{proof}
  Let $E,A,B,H$ be a set of real matrices whose fiber has
  dimension equal to $d-n-1$. (Since the set of real matrices
  inside $\C^{2m(m+n)}$ is not contained in an algebraic set,
  such real matrices exist.) Let $\cZ$ be the induced center.  If
  the degree of $\bar{\cL}$ is odd then the finite morphism
  $\pi_1:\bar{\cL}\rightarrow \pn{d}$ is surjective over the real
  numbers. Indeed over the complex numbers the inverse image
  $\pi_1^{-1}(y)\subset \bar{\cL}$ represents a finite set of
  complex conjugate points for every real point $y\in\pn{d}$.
  But then also $\pi_2$ and $\chi$ are surjective over the reals.
\end{proof}

\Section{The degree of  $\bar{\mathcal{L}}$ in some special
  situations} \label{Sec-special}

In this and in the next section we will assume that $\K$ is an
algebraically closed field of characteristic zero. We will show
in a moment that the compactification $\bar{\cL}$ is in many
cases isomorphic to the product of some Schubert varieties. This
will allow us to compute the degree of $\bar{\cL}\subset\pn{N}$
in these cases.

For the convenience of the reader we summarize the basic notions.
More details can be found in~\cite{ho52,ro97} and~\cite[Chapter
14]{fu84}.

Consider a flag
$$
{\mathcal F}: \ \{ 0\}\subset V_{1}\subset V_{2} \subset\ldots
\subset V_{m+n}= {\K}^{m+n}
$$
where we assume that $\dim V_{q}=q$ for $ q=1,\dots ,m+n$.
Let $\nu= ( \nu_{1},\ldots ,\nu_{m})$ be an ordered index set
satisfying
$$
1\leq \nu_{1} <\ldots <\nu_{m} \leq m+n.
$$
With respect to the flag ${\mathcal F}$ one defines the
Schubert variety
$$
S(\nu_1,\ldots,\nu_m) := \{ W\in \G(m,\K^{m+n})\mid
\dim(W\bigcap V_{\nu_{k}}) \geq k \ \mbox{ for } k=1,\ldots ,m\}
$$
and the Schubert cell
$$
C(\nu_1,\ldots,\nu_m):=\{W\in S(\nu_1,\ldots,\nu_m)\mid \dim
(W\bigcap V_{\nu_{k}-1})=k-1 ; \mbox{ for } k=1,\ldots ,m\}.
$$
The closure of the Schubert cell $C(\nu_1,\ldots,\nu_m)$
inside the variety $\G(m,\K^{m+n})\subset \pn{N}$ is equal to the
Schubert variety $S(\nu_1,\ldots,\nu_m)$. By definition,
$S(\nu_1,\ldots,\nu_m)$ is a projective variety.  There is a well
known formula for the degree of a Schubert variety~\cite[Chapter
XIV, \S 6, (7)]{ho52}:
$$
\deg S(\nu_1,\dots,\nu_k)=(\sum_i(\nu_i-i))!
\frac{\prod_{j>i}(\nu_j-\nu_i)}{\prod_i(\nu_i-1)!}.
$$

Let ${\mathcal B}:=\{ v_1,\ldots,v_{m+n}\} \subset \K^{m+n}$ be a
basis which is compatible with the flag ${\mathcal F}$.  In other
words this basis has the property that $V_i={\rm
  span}(v_1,\dots,v_i)$. With respect to the basis ${\mathcal B}$
one can represent the Schubert cell $ C(\nu_1,\ldots,\nu_m)$ as
the set of all $m$-dimensional subspaces in $\K^{m+n}$ which are
the rowspaces of a matrix of the form:
\begin{equation}
  \label{cell}
 \left[\begin{array}{cccccccccccccccc} \ast &\cdots&\ast&1 &0
      &\cdots&0 &0
      &\cdots&0   &\cdots&0   &0&0&\cdots&0\\
      \ast &\cdots&\ast&0 &\ast&\cdots&\ast &1
      &\cdots&0   &\cdots&0   &0&0&\cdots&0\\
      \vdots& &\vdots&\vdots&\vdots& &\vdots&\vdots&
      &\vdots&    &\vdots&\vdots&\vdots&&\vdots\\
      \ast &\cdots&\ast&0 &\ast&\cdots&\ast
      &0&\cdots&\ast&\cdots&\ast&1 &0&\cdots&0
  \end{array}\right]
\end{equation}
where the $1$'s are in the columns $\nu_1,\dots,\nu_m$.

The cell $C(\nu_1,\ldots,\nu_m)$ is isomorphic to $\K^d$, where
$d=\sum_{i=0}^m(\nu_i-i)$. In particular the cell
$C(\nu_1,\ldots,\nu_m)$ is isomorphic to every subspace
$\cL\subset\matmn$ having dimension $\dim\cL=d$. In general it is
not true that the closures $S(\nu_1,\ldots,\nu_m)\subset\pn{N}$
and $\bar{\cL}\subset\pn{N}$ are isomorphic. This happens however
in the following situation:

Let $E_{i,j}$ be the $m\times n$ matrix whose $i,j$-entry is $1$
and all the other entries are $0$. Let $\mu= ( \mu_{1},\ldots
,\mu_{m})$ be an ordered index sets satisfying
$$
0\leq \mu_{1} \leq\ldots \leq\mu_{m} \leq n.
$$
\begin{defn}
  $\cL\subset\matmn$ is called a {\em lower left filled linear
    space} of type $\mu$ if $\cL$ is spanned by the matrices
  $$
  E_{i,j} \,\,\, \mbox{ for }\, \ j\leq\mu_i,\ i=1,\ldots,m.
  $$
\end{defn}
\begin{lem}
  If $\cL\subset\matmn$ is a lower left filled linear space of
  type $\mu$ then $\bar{\cL}$ is isomorphic to the Schubert
  variety $S(\mu_{1}+1,\mu_{2}+2,\ldots,\mu_{m}+m)$.
\end{lem}
\begin{proof}
  Let $\nu_i:=\mu_i+1$,\ $i=1,\ldots,m$. There is a fixed
  $(m+n)\times (m+n)$ permutation matrix $P$ such that the set
  $$
  \left\{ \left[ I_m | F\right] P\mid F\in\cL\right\}\subset%
  M\! at_{m\times (m+n)}
  $$
  is equal to the cell $C(\nu_1,\ldots,\nu_m)$ described
  in\eqr{cell}. The linear transformation $P\in Gl_{m+n}$ extends
  to a linear transformation in ${\mathbb P}\left(\wedge^{m}
    \K^{m+n}\right)=\pn{N}$ and this linear transformation maps
  $\bar{\cL}$ isomorphically onto $S(\nu_1,\ldots,\nu_m)$.
\end{proof}

The proof of the lemma shows in particular that permutations of
the columns inside $\matmn$ result in isomorphic
compactifications. The following lemma shows that a broader range
of transformations do not change the topological properties of
the compactification.

\begin{lem}                \label{lem33}
  Assume there are subspaces $\cL_1,\cL_2\subset\mat$. If there
  are linear transformations $S\in Gl_m$ and $T\in Gl_n$ such
  that $\cL_2=S\cL_1 T^{-1}$ then there exists an automorphism of
  $\pn{N}$ which maps the compactification $\overline{\cL_1}$
  isomorphically onto the compactification $\overline{\cL_2}$.
\end{lem}
\begin{proof}
  $$
  \left[ I_m\mid S\cL_1 T^{-1}\right]=S\left[ I_m\mid \cL_1
  \right]\vier{S^{-1}}{0}{0}{T^{-1}}.
  $$
  The matrix to the right, an element of $GL_{m+n}$, induces a
  linear transformation on the projective space ${\mathbb
    P}\left(\wedge^{m} \K^{m+n}\right)=\pn{N}$ which maps
  $\overline{\cL_1}$ onto $\overline{\cL_2}$.
\end{proof}

\begin{thm}                  \label{prod}
  Assume there are linear transformations $S\in Gl_m$ and $T\in Gl_n$
  such that
  $$
  S\cL T^{-1}= \left(
\begin{array}{ccc}
\cL_1 & 0 & 0\\
0 & \ddots & 0\\
0&       0 &\cL_k
\end{array}
\right),
$$
where each $\cL_l,l=1,\ldots,k$ is the space of $m_l\times
n_l$ lower left filled matrices of type $\mu^{l}$:
$$
0\leq \mu_{1}^{l}\leq\cdots\leq \mu_{m_l}^{l}\leq n_l.
$$
Then $\bar{\cL}$ is isomorphic to the product of Schubert
varieties
$$
S(\mu_{1}^{1}+1,\mu_{2}^1+2,\dots,\mu_{m_1}^1+m_1)\times\cdots
\times S(\mu_{1}^k+1,\mu_{2}^k+2,\dots,\mu_{m_k}^k+m_k)
$$
and
$$
\deg\bar{\cL}=\left(\sum_{i,l}\mu_{i}^l\right)!\frac{\displaystyle{
    \prod_{i,l_r >l_s}\left(\mu_{i}^{l_r}+l_r
      -\mu_{i}^{l_s}-l_s\right)}}{\displaystyle{\prod_{i,l}
    \left(\mu_{i}^l+l-1\right)!}}.
$$
\end{thm}

\begin{proof}
  The closure of $[\cL_l,I_{m_l}]$ in the Grassmannian variety
  $\G (m_l,\K^{m_l+n_l})$ is the Schubert variety
  $S(\mu_1^l+1,\dots,\mu_{m_l}^l+m_l)$, and $\bar{\cL}$ is a
  product of Schubert varieties.
  
  The degree formula of a product of projective varieties under
  the Segre embedding~\cite[Proposition~2.1]{wa94a} is given by
  $$
  \deg Z_1\times \cdots\times Z_k=\frac{(\sum_i \dim
    Z_i)!}{\prod_i (\dim Z_i)!}\prod_i \deg Z_i.
  $$
  Combining these formulas gives the result.
\end{proof}

\begin{cor}           \label{BrBy2}
  When $\mu_1^1=\cdots=\mu_{m_1}^1=n_1$ and $\mu_i^l=0$ for
  $l>1$, then the compactification $\bar{\cL}=\G
  (m_1,\K^{m_1+n_1})$ and its degree is
  $$
  \frac{(m_1n_1)! 1!2!\cdots (m_1-1)!}{n_1!(n_1+1)!\cdots
    (n_1+m_1-1)!}.
  $$
\end{cor}

Using Lemma~\ref{lem33} and Corollary~\ref{BrBy2} we can deduce
Theorem~\ref{BrBy}, the result of Brockett and Byrnes. For this
assume that $\cL=\{ BFC\mid F\in\matmp\}$. Without loss of
generality we can assume that $B,C$ have full rank, $\rank B=m$
and $\rank C=p$. (Theorem~\ref{BrBy} assumes genericity!) There
are invertible matrices $S,T$ such that $SB=\zwei{I_m}{0}$ and
$CT^{-1}=\left[ I_p|0\right]$. It follows that
$$
S\cL T^{-1}=\left\{ \vier{F}{0}{0}{0}\in\mat\mid
  F\in\matmp\right\} .
$$
According to Lemma~\ref{lem33} and Corollary~\ref{BrBy2} the
compactification is isomorphic to the Grassmannian $\G
(m,\K^{m+p})$ as predicted by Theorem~\ref{BrBy}. In order to
fully prove Theorem~\ref{BrBy} it remains to be shown that for a
generic set of matrices $A\in\mat$ the system is
$\cL$-nondegenerate as soon as $n=mp$.

\begin{cor}           \label{friedl}
  When $m_l=n_l=1$ and $\mu_1^l=1$ for all $l$, then
  $\bar{\cL}=\prod_{i=1}^n\pn{1}=\pn{1}\times\cdots\times \pn{1}$
  and its degree is
  $$
  n!.
  $$
\end{cor}

Corollary~\ref{friedl} covers a result first studied by
Friedland~\cite{fr72,fr77a}. Indeed the subspace $\cL\subset\mat$
corresponds in this case exactly to the set of diagonal matrices.
By Theorem~\ref{main} we know that for a generic set of matrices
$E,A,B,H$ the characteristic map $\chi$ is a finite morphism of
mapping degree $m!$.  Friedland~\cite{fr72,fr77a} and Byrnes and
Wang~\cite{by93} did show that the set of all matrices of the
form $I_n,A,I_n,0$ belongs to this generic set. We therefore have
the result:
\begin{thm}[\cite{by93,fr72,fr77a}]                \label{friedlThm}
  Let $\cL\subset\mat$ be the set of all diagonal matrices
  defined over an algebraically closed field $\K$. If $A\in\mat$
  is an arbitrary matrix and $\varphi\in\K[s]$ is an arbitrary
  monic polynomial of degree $n$ then there are exactly $n!$
  diagonal matrices $F\in\cL$ (when counted with multiplicity)
  such that the matrix $A+F$ has characteristic
  polynomial~$\varphi(s)$.
\end{thm}

\Section{The degree of $\bar{\mathcal{L}}$
  in the generic situation} \label{Sec-generic}

In the previous section we computed the degree of the variety
$\bar{\cL}$ in many special cases. The set of all subspace
$\cL\subset\matmn$ having the property that $\dim\cL=d$ can be
identified with the Grassmannian variety $\G(d,\K^{mn})$. The
degree attains its maximal value on a Zariski open subset of
$\G(d,\K^{mn})$. This largest possible degree is sometimes
referred to as the generic degree. In this section we will
determine this generic degree if $d=m=n$. The result is as
follows:

\begin{thm}\label{thm4.1}
  Let $\K$ be an algebraically closed field of characteristic
  zero.  There is a generic subset $U\subset \G(n,\K^{n^2})$ such
  that the compactification $\bar{\cL}\subset {\mathbb
    P}\left(\wedge^{n} \K^{2n}\right)={\mathbb P}^N$ of every
  element $\cL\in U$ has degree $n(n-1)^{n-1}.$ This is also
  equal to the maximal possible degree among all varieties
  $\bar{\cL}$ with $\dim\cL=n$ and $\cL\subset\mat$.
\end{thm}

The proof of this theorem will require a fair amount of algebraic
geometry. In particular, our proof involves a ``blowing-up
method.'' This is an important method in algebraic geometry and
is the main tool in the resolution of singularities of algebraic
varieties, or in the elimination of the points of indeterminacy
of  rational maps in consideration.  (A rational map is a
morphism which is only defined on some open and (Zariski-) dense
subset.) This
blowing-up construction enables us to compute the degree of the
rational map. The interested reader may want to
consult~\cite{gr78,ha77} for the notation and basic facts on
blowing-up.

In order to make the proof more understandable we will explain it
first in the specific examples $n=3$ and $n=4$. Thereafter we
will give the general proof.

\begin{exmp}\label{exn3}
  Let
  $$
  \cL= \left[
\begin{array}{ccc}
z_1&z_2&z_3\\
0&z_1&z_2\\
z_3&0&z_1
\end{array}
\right],
$$
be a $3$-dimensional linear subspace in $M\! at_{3\times 3}.$
The full size minors of $[I,\cL]$ give the following 20
coordinates:
$$
(1,z_1,z_2,z_3,z_1^2,z_2 z_3,-z_1 z_3, z_3^2-z_1^2, z_2^2-z_1
z_3,\ldots ,z_1^3+z_2^2 z_3-z_1 z_3^2).
$$
By adding another variable $z_0$ to compactify $\cL$ and
homogenize the coordinates, we get
$$
(z_0^3,z_0^2 z_1,z_0^2 z_2,z_0^2 z_3,z_0 z_1^2,z_0 z_2z_3,
\ldots ,z_1^3+z_2^2 z_3-z_1 z_3^2).
$$
Let
$$
\phi : \pn{3} \longrightarrow \pn{19}
$$
be the rational map defined by the above set of degree 3 homogeneous
polynomials on $\pn{3}$, say ${\mathcal D},$ which is a sublinear
system of the complete linear system $| \so_{\pn{3}} (3)|$.
In general, 
$| \so_{\pn{n}} (d)|$ determines a morphism, which is called
the $d$-uple embedding, from $\pn{n}$ to $\pn{N}$ where
$N= {{n+d}\choose n}-1$, defined
by the algebraically independent  homogeneous polynomials of degree $d$
in $n+1$ variables.

Note that $\phi$ is not defined on the cubic curve $B:=\{z\in
\pn{3}|z_0=0, z_1^3+z_2^2 z_3-z_1 z_3^2=0\}$ in $\pn{3}$. The
curve $B$ coincides with the indeterminacy locus of the rational
map $\phi$ (also scheme theoretically to be precise).  We denote
$\phi$ by $|{\mathcal D} -B|$.  Note also that $B$ is nonsingular
and irreducible. Let
$$
\pi : \widetilde{\pn{3}} \rightarrow \pn{3}
$$
be the blowing-up of $I_B$, the ideal of $B$ in $\pn{3}$.  
Since we blew-up the smooth curve $B$, $\widetilde{\pn{3}}$ is a
projective manifold containing the smooth exceptional divisor
$\tB:=\pi ^{-1} (B)$. Let
$\tD:=\pi^* {\mathcal D}$ be the pulled-back sublinear system on
$\widetilde{\pn{3}}$.
$\tB$ is isomorphic to
$\pn{}(N_{\!B\!/{\pn{3}}}^\vee )$, the projective space bundle of
hyperplanes in the conormal bundle $N_{\!B\!/{\pn{3}}}^\vee$
of rank equal to 2, codimension of $B$ in $\pn{3}$. 
There is the natural
 projection morphism $p:\pn{}(N_{\!B\!/{\pn{3}}}^\vee
)\rightarrow B$, and  the tautological line bundle $\xi$ on
$\pn{}(N_{\!B\!/{\pn{3}}}^\vee )$ whose restriction to the fiber
$\pn{1}$ of $p$ is isomorphic to $\so_{\pn{1}}(1)$.
By the argument  as in \cite[II, Example 7.17.3]{ha77}, we have
a well-defined morphism $\widetilde{\phi}=|\tD -\tB|:
\widetilde{\pn{3}}\rightarrow \pn{19}$ factoring through
$\phi$:
$$
\xymatrix{
  \widetilde{\pn{3}}\ar[d]_{\pi}\ar[dr]^{\widetilde{\phi}} \\
  \pn{3}\ar[r]^{\phi}& \pn{19}.}
$$
Since $\cL$ can be identified as $\{z=(z_0, z_1, z_2, z_3)\in
\pn{3}|z_0\ne 0\}$, $\cL$ lies in $\pn{3}-B$ on which $\pi$ is an
isomorphism. So the degree of the closure $\overline{\phi(\cL )}$
of ${\phi(\cL )}$ is equal to the 
degree of 
${\widetilde{\phi}} (\widetilde{\pn{3}})$ which is equal to the
self-intersection number
$$
(\tD- \tB )^3 =\tD ^3 - 3 \tD ^2 \cdot \tB +3\tD \cdot \tB^2
-\tB^3.
$$
\begin{itemize}
\item[(i)] $\tD ^3=3^3$ since ${\mathcal D}^3$ is the
  intersection number of three hypersurfaces of degree 3 on
  $\pn{3}$ which are generic elements in ${\mathcal D}$.
\item[(ii)] $\tD ^2 \cdot \tB=0$ since ${\mathcal D}^2 $ is a
  curve in $\pn{3}$ and generic enough not to meet the cubic
  curve $B$.
  
  Let $E$ be a vector bundle of rank $r$ on a nonsingular
  projective variety $X$.  In the next computation, we use the
  following `Chern relation' in the cohomology ring of the
  projective space bundle $p: \pn{} (E)\rightarrow X$ which gives
  a relation between the Chern classes of $E$ and the tautological
  line bundle $\xi$ on $\pn{} (E)$:
  $$
  \sum_{i=0}^r (-1)^i p^* c_i (E) \xi^{r-i}=0
  $$
  in $H^{2r} (\pn{} (E))$.  We will suppress the pull-back
  $p^*$.
\item[(iii)] $\tB^3 = \tB|_{\tB} \cdot \tB|_{\tB} = {(-\xi )}^2$.
  By the Chern relation for $\pn{}(N_{\!B\!/{\pn{3}}}^\vee)$
  $$
  \xi^2 - c_1 (N_{\!B\!/{\pn{3}}}^\vee)\,\xi +
  c_2(N_{\!B\!/{\pn{3}}}^\vee )=0,
  $$
  and by the fact that $c_2(N_{\!B\!/{\pn{3}}}^\vee)$
  restricted to $B$ is  automatically trivial since $\dim B <2$, 
  we get $\xi^2 = c_1
  (N_{\!B\!/{\pn{3}}}^\vee )\,\xi=c_1 (N_{\!B\!/{\pn{3}}}^\vee
  )$.  Now observe that $N_{\!B\!/{\pn{3}}}\cong\so_{\pn{3}} (1)
  \oplus \so_{\pn{3}} (3)$ since $B$ is a complete intersection
  of two hypersurfaces of degree 1 and 3.  Therefore
  $$
  \xi^2=\deg\, (\wedge^2 N_{\!B\!/{\pn{3}}}^\vee)|_B = -\deg\,
  \so _{B} (4) =-12,
  $$
  which yields $\tB^3=-12$.
\item[(iv)] Finally,
\begin{eqnarray*}
\tD \cdot \tB^2 & = &\tD|_{\tB} \cdot\tB|_{\tB} \\
& = & (9\, \mbox{- fibers of }\pi|_{\tB} ) \cdot (-\xi) \\
& = &-9
\end{eqnarray*}
\end{itemize}
Therefore $(\tD- \tB )^3=12=3\cdot 2^2$.
\end{exmp}

\begin{exmp} Let $\cL=\{M(z_1,\dots,z_4)\}\subset M\! at_{4\times 4}$ be a
  $4$-dimensional subspace where $M(z_1,\dots,z_4)$ represents a
  $4\times 4$ matrix whose entries are given by linear
  polynomials in $z_1,\dots z_4$.  Then the indeterminacy locus
  $$
  B=\{z=(z_0, z_1, z_2, z_3, z_4)\in \pn{4}|z_0=0, f(z_1,\dots,z_4)=0) \}
  $$
  where $f(z_1,\dots,z_4)=\det M(z_1,\dots,z_4)$ is a
  homogeneous polynomial of degree $4$.  We may assume that $B$
  is nonsingular and irreducible for the moment; the reason for
  it will be given in the proof of Theorem~\ref{thm4.1}.  Such an
  assumption is necessary in order to carry out the computation
  using the normal bundle $N_{\!B\!/{\pn{4}}}$ of $B$ in $\pn{4}$
  and to use the Chern relations as in the previous case.  Let
  $\phi : \pn{4} \longrightarrow \pn{{8 \choose 4}-1}$ be the
  rational map defined by the sublinear system $\cal D$ of the
  complete linear system $|\so_{\pn{4}} (4)|$, which is obtained
  by taking the Pl\"ucker coordinates for $\cL$.  Now we have the
  following diagram:
  $$
  \xymatrix{
  \widetilde{\pn{4}}\ar[d]_{\pi}\ar[dr]^{\widetilde{\phi}} \\
  \pn{4}\ar[r]^{\phi}& \pn{{8 \choose 4}}}
  $$
  where $\pi$ is the blowing-up of $I_B$, $\phi= |{\cal D} -B|$
  and $\widetilde{\phi}= |\tD -\tB |$ as in the case $n=3$.  Now we
  want to calculate the self-intersection number
  $$
  (\tD- \tB )^4 =\tD ^4 - {4 \choose 1} \tD ^{3} \cdot \tB + {4
  \choose 2}\tD^{2} \cdot \tB^2 -{4 \choose 3} \tD \cdot \tB^3 +
  \tB^4.
  $$
  \begin{itemize}
  \item[(i)] $\tD^4=4^4$.

  \item[(ii)] $\tD^3 \cdot \tB =0$ since $\tD ^3 $ is a curve in
  $\pn{4}$ and generic enough not to meet the quadric surface
  $B$.
  
\item[(iii)] $\tD \cdot \tB^3=4\,\pi^* H \cdot \tB^3$ where $H\in
  |\so_{\pn{4}} (1)|$.  $\pi^* H \cdot \tB^3=\tB^3 |_{\pi^* H}$,
  i.e. $\tB^3$ restricted to the pull-back of the curve $C:=H\cap
  B$ of degree $4$ in $H\cong\pn{3}$.  Let us still denote it by
  $\tB^3$.  Note that by Bertini's theorem \cite[II, Theorem
  8.18]{ha77}, for generic $H$, $C$ is nonsingular and
  irreducible.  Then on $\pi^{-1}
  (C)\cong\pn{}(N_{\!C\!/{\pn{3}}}^\vee)$,
  $$
  \tB^3=\tB |_{\tB} \cdot \tB |_{\tB}=(-\eta)^2
  $$
  where $\eta$ is the tautological line bundle on
  $\pn{}(N_{\!C\!/{\pn{3}}}^\vee)$. Observe that
  $N_{\!C\!/{\pn{3}}}\cong\so_{\pn{3}}(1)\oplus \so_{\pn{3}}(4)$.
  Now we use the Chern relation for
  $\pn{}(N_{\!C\!/{\pn{3}}}^\vee)$: $\eta^2 - c_1
  (N_{\!C=\!/{\pn{3}}}^\vee)\eta=0$ as in the previous case.
  Therefore $\eta^2=c_1(N_{\!C\!/{\pn{3}}}^\vee)=-5\cdot 4=-20$
  and $\tD \cdot \tB ^3=- 4^2\cdot 5$.

\item[(iv)] $\tD^2 \cdot\tB^2=16 \,\pi^* H^2|_{\tB} .
  \tB|_{\tB}=-4^3$ as in the previous case.

\item[(v)] Finally, $\tB^4= (\tB|_{\tB})^3=(-\xi)^3$ where $\xi$
  is the tautological line bundle on the exceptional divisor
  $\tB\cong\pn{}(N_{\!B\!/{\pn{4}}}^\vee)$.  By the Chern
  relation
  $$
  \xi^3 - c_1 (N_{\!B\!/{\pn{4}}}^\vee )\,\xi^2 +
  c_2(N_{\!B\!/{\pn{4}}}^\vee )\,\xi=0,
  $$
  and by the fact that $N_{\!B\!/{\pn{4}}}\cong\so_{\pn{4}}
  (4) \oplus \so_{\pn{4}} (1)$, we have
\begin{eqnarray*}
c_1 (N_{\!B\!/{\pn{4}}}^\vee ) \xi^2
   &= & (\pi^* \so_{\pn{4}}(-5)) \xi^2\\
   &= &-5 (\pi^* \so_{\pn{4}}(1)) (-\tB|_{\tB} )^2\\
   &= &(-5) \cdot \eta^2 \\
   &= & (-5) \cdot (-20)\\
   &= & 5^2\cdot 4,
\end{eqnarray*}
and
$$
c_2(N_{\!B\!/{\pn{4}}}^\vee )\, \xi = 4 \,(\pi^* \so_{\pn{4}}
(1))^2\, \xi =(4^2\,\mbox{- fibers of }\pi|_{\tB} )\cdot\xi =4^2,
$$
which yields $\tB^4= -4(4^2 +4+1)$.
\end{itemize}
Therefore
$$
\tB^4=(\tD- \tB )^4 =4\cdot 3^3.
$$
\end{exmp}

\begin{proof}[Proof of Theorem~\ref{thm4.1}]
  Let $\cL=\{M(z_1,\dots,z_n)\}\subset \mat$ be an
  $n$-dimensional subspace where $M(z_1,\dots,z_n)$ represents a
  $n\times n$ matrix whose entries are given by linear
  polynomials in $z_1,\dots z_n$.  The Pl\"{u}cker coordinates of
  $[I,M(z_1,\dots,z_n)]$ define a polynomial map from $\A_{\K}^n$
  to $\A_{\K}^{{2n \choose n}-1}\subset \pn{{2n \choose n}-1}$.
  Homogenizing the map, we have a map $\phi : \pn{n}
  \longrightarrow \pn{{2n \choose n}-1}$ defined by the sublinear
  system $\cal D$ of $|\so_{\pn{n}} (n)|$.  The restriction of
  $\phi$ to $\cL$ is the Pl\"ucker embedding which is described
  in the beginning of Section 2.  However, the map $\phi$ is not
  well-defined on a subvariety of codimension 2, i.e.  $\phi$ is
  not a ``morphism," but a ``rational map" with non-empty
  indeterminacy locus.  To eliminate this indeterminacy locus we
  will construct a ``blow-up."
  
  Let
  $$
  B=\{z=(z_0, z_1, \dots, z_n)\in \pn{n}|z_0=0, f(z_1,\dots,z_n)=0) \}
  $$
  be the indeterminacy locus of pure codimension 2 where
  $f(z_1,\dots,z_n)=\det M(z_1,\dots,z_n)$ is homogeneous of
  degree $n$. We may assume that $B$ is
  nonsingular and irreducible.  If not, consider
  $$
  p: \pn{n}\times \A_{\K}^1\rightarrow \A_{\K}^1,
  $$
  the trivial flat family of $\pn{n}$ over $\A_{\K}^1$ (or the
  projection onto $\A_{\K}^1$) where $\A_{\K}^1$ is the affine
  line over $\K$.  Let $\pn{n}_t\cong\pn{n}$ be the fiber over
  $t\in\A_{\K}^1$, and consider the subvariety $B_t$ of
  $\pn{n}_t$ given by $(z_0=0, f+tg=0)$ where $g\in \K
  [z_1,\ldots , z_n]$ is homogeneous of degree $n$.  Such $\{
  B_t\}$ forms a flat family over $\A_{\K}^1$ by varying $t\in
  \A_{\K}^1$ with $B_0=B$.  Moreover, we can choose an
  appropriate $g$ so that the generic $B_t$ is nonsingular and
  irreducible.  Let $I_{B_t}$ be the ideal of $B_t$ in $\pn{n}_t$
  for $t\in \A_{\K}^1$.  Let ${\mathcal D_t}$ be the sublinear
  system of $|\so_{\pn{n}}(n)|$ consisting of the same polynomials
  of degree $n$ as for ${\mathcal D}$, except that $f$ is replaced by
  $f+tg.$ For each $t\in \A_{\K}^1$, we take the blowing-up
  $\pi_t$ of $I_{B_t}$ with the exceptional divisor
  $\widetilde{B_t}$, satisfying the following diagram:
  $$
  \xymatrix{
  \widetilde{\pn{n}_t}\ar[d]_{\pi_t}\ar[dr]^{\widetilde{\phi_t}} \\
  \pn{n}_t\ar[r]^{\phi_t}& \pn{{2n  \choose n}-1}}
  $$
  with $\phi_t=|{\mathcal D_t}-B_t|$ and
  $\widetilde{\phi_t}=|\widetilde{{\mathcal
      D}_t}-\widetilde{B_t}|$ where 
      ${\widetilde{\mathcal D}_t}:=\pi_t^* {\mathcal D}_t$
      is the pulled-back sublinear system on $\widetilde{\pn{n}_t}.$
  On the other hand, we take the blowing-up of $I$, the ideal of
  the subscheme $\{ B_t\}_{t\in\A_{\K}^1}$ in $\pn{n}\times
  \A_{\K}^1$:
$$
\pi :\widetilde{\pn{n}\times\A_{\K}^1}\rightarrow \pn{n}\times
\A_{\K}^1.
$$
with the exceptional divisor $E$.  Since $\A_{\K}^1$ is
a smooth curve, the composition
$$
p\circ\pi:
\widetilde{\pn{n}\times\A_{\K}^1}\rightarrow\A_{\K}^1
$$
is flat (cf. \cite[Appendix B.6.7]{fu84}).  The fiber of
$p\circ\pi$ over $t=0$ is $\widetilde{\pn{n}}$, i.e. the blowing
up of $I_0$, the ideal of $B$ in $\pn{n}$. Note that the
restriction of the exceptional divisor $E$ to a fiber,
$E|_{(p\circ\pi)^{-1}}(t)$, is $\widetilde{B_t}$.  The flatness
of $p\circ\pi$ assures that certain numerical invariants remain
constant in the family (see \cite[III, Theorem 9.9]{ha77}).  In
particular, the self-intersection
$(\widetilde{D_t}-\widetilde{B_t})^n$ is independent of $t\in \A_{\K}^1$.
This is why  we can assume that $B$ is  irreducible and nonsingular.

Let us proceed with the proof of the theorem.  We write
$\phi=|{\cal D} -B|$.  We have the following diagram:
$$
  \xymatrix{
  \widetilde{\pn{n}}\ar[d]_{\pi}\ar[dr]^{\widetilde{\phi}} \\
  \pn{n}\ar[r]^{\phi}& \pn{{2n  \choose n}-1}}
$$
where $\pi$ is the blowing-up of $I_B$.  $\so_{\pn{n}} (1)
\oplus \so_{\pn{n}} (n)$.  By using the Chern relation
$$
\xi^{n-1}-c_1(N_{\!B\!/{\pn{n}}}^\vee )\,\xi^{n-2}+\cdots
+(-1)^{n-2} c_{n-2}(N_{\!B\!/{\pn{n}}}^\vee )\,\xi=0
$$
inductively as in the previous cases, we get the following:
$$
\begin{array}{lcl}
  \tD^n              &= &n^n\\
  \tD ^{n-1} \cdot \tB   &= &0\\
  \tD ^{n-2} \cdot \tB^2 &= &-n^{n-1}\\
  \tD ^{n-3} \cdot \tB^3 &= &-n^{n-2} (n+1)\\
  \tD ^{n-4} \cdot \tB^4 &= &-n^{n-3} (n^2 + n+1)\\
  \tD ^{n-5} \cdot \tB^5 &= &-n^{n-4} (n^3 + n^2 +n+1)\\
  &\vdots &\\
  \tD ^{0} \cdot \tB^n &= &- n(n^{n-2} + n^{n-3} +\cdots
  +n+1),
\end{array}
$$
which yields
\begin{eqnarray*}
(\tD-\tB)^n &=&\sum_{k=0}^n(-1)^{k}{n\choose k}\tD^{n-k}\cdot\tB^k\\
&= & n^n +\sum_{k=2}^n \left[ (-1)^{k+1}{n\choose k} n^{n-k+1}\big(
\sum_{i=0}^{k-2} n^i \big)\right] \\
&= &n(n-1)^{n-1}.
\end{eqnarray*}
\end{proof}

Theorem~\ref{thm4.1} will allow us to answer
Problem~\ref{problem} in the ``generic situation'' if $m=n=d$:

\begin{thm}
  Assume $\K$ is an algebraically closed field of characteristic
  zero. Let $\cL\subset\mat$ be a ``generic subspace'', let
  $E,A,B,H$ be a ``generic set of matrices'' and let
  $\varphi(s)\in\K^n$ be a ``generic monic polynomial'' of degree
  $n$. Then there exist exactly $n(n-1)^{n-1}$ different feedback
  laws $F\in\cL$ such that $\det \left[
    s(E+HF)+(A+BF)\right]=\varphi(s)$. In particular
 the system\eqr{dis-sys} is generically pole assignable in the class
  of feedback compensators~$\cL$.
\end{thm}
\begin{proof}
  Consider the characteristic map $\chi$ introduced
  in\eqr{central}. According to Theorem~\ref{main} $\chi$ is a
  finite morphism of degree $\bar{\cL}$. By Theorem~\ref{thm4.1}
  the degree of $\bar{\cL}$ is $n(n-1)^{n-1}$. For a generic set
  of polynomials $\varphi(s)\in\K^n\subset\pn{n}$ the inverse
  image $\chi^{-1}(\varphi(s))$ contains
  $\deg\left(\bar{\cL}\right)$ different solutions and all these
  solutions are contained in $\cL\subset\bar{\cL}$.
\end{proof}


\end{document}